\documentclass[10pt]{amsart}
\usepackage{amsmath,amssymb,amscd,enumerate,bbm}
\usepackage{cite}

\renewcommand{\MR}[1]{}

\newcommand{\titl}{A REMARK ON MINIMAL FANO THREEFOLDS}

\title{{\titl}}
\author{V. Golyshev}
\date{}

\newcommand{\cal}{\mathcal}

\def\C{{\Bbb{C}}}
\def\Q{{\Bbb{Q}}}
\def\Z{{\Bbb{Z}}}

\def\L{{\cal{L}}}

\def\OO{{\mathcal{O}}}
\renewcommand{\phi}{{\varphi}}

\def\Hom{{\rm H}{\rm o}{\rm m}}

\newcommand{\cD}{{\cal D}}

\newcommand{\Auteq}{{AutEq}}
\newcommand{\Ob}{{Ob}}
\newcommand{\Ker}{\mathop{Ker}}

\newcommand{\dimens}{{\mathrm{dim \;} }}

\newcommand{\mapto}{\mapsto}

\newcounter{pphcounter}[section]
\renewcommand{\thepphcounter}{\thesection.\arabic{pphcounter}}
\newcommand{\pph}{\bigskip \refstepcounter{pphcounter}
    \bf  \thepphcounter. \rm}

\makeatletter
\@addtoreset{equation}{pphcounter}
\makeatother

\newcommand{\defi}{\bf Definition. \rm}

\newcommand{\propo}{\bf Proposition. \rm}
\newcommand{\theo}{\bf Theorem. \rm}

\renewcommand{\proof}{\bigskip \bf Proof. \rm}

\def\A1{{{\Bbb{A}}^1}}
\def\P1{{{\Bbb{P}}^1}}

\def\Gm{{\bf G_m}}

\def\SL2{{\mathrm SL2}}

\renewcommand{\P}{{\Bbb{P}}}

\newcommand{\F}{{\cal F}}

\newcommand{\bdx}{{\mathbf x}}
\newcommand{\bdE}{{\mathbf E}}

\begin{document}

\begin{center}
\maketitle

\end{center}

\bigskip

\parbox{340pt}{\small \bf Abstract. \rm
We prove in the case of minimal Fano
 threefolds a conjecture stated by Dubrovin at the ICM 1998 in Berlin.
The conjecture predicts that the symmetrized/alternated Euler
characteristic pairing on $K_0$ of a Fano variety with an exceptional
collection expressed in the basis of the classes of the exceptional
objects coincides with the intersection pairing of the vanishing cycles
in Dubrovin's second connection. We show that the conjecture holds for
$V_{22}$, a minimal Fano threefold of anticanonical degree~$22$, and
for $V_5$, the minimal Fano threefold of anticanonical degree~$40$, by applying
the modularity result for the rank $1$
Fano threefolds  established in \cite{Go-CPMD}. The
truth of the conjecture for $\P ^3$ and
the three--dimensional quadric is known; we consider these cases
for the sake of completeness.}

\bigskip
\bigskip
\bigskip

\section{The
conjecture}

\pph {\bf Exceptional collections and autoequivalences.} Let $F$ be an
algebraic variety. By   $\cD^b(F)$ denote the bounded derived
category of coherent sheaves on $F$.

Any derived category  $\cD$ is triangulated. This means that
given are the shift functor $[1] : \cD \to \cD$ which is an additive
autoequivalence, and the class of distinguished triangles
$
X \to Y \to Z \to X[1],
$
which satisfies the standard axioms.
An additive functor  $\F :\cD \to \cD'$ between two triangulated
categories is said to be {\it exact} if it commutes with the
shift functor and sends distinguished triangles to distinguished
triangles. The group of isomorphism classes of
exact equivalences  $\F : \cD^b(F) \to \cD^b(F)$ is called {\it the
group of autoequivalences of $\cD^b(F)$\/} and denoted
$\Auteq(\cD^b(F))$.

By  $\Hom^i(X,Y)$ denote $\Hom(X,Y[i])$.  An object
$E \in \mathrm{Ob}(\cD^b(F))$ is said to be {\it exceptional} if it
satisfies
$$
\Hom^i(E,E)=0\quad\text{when}\quad i>0,\quad \Hom^0(E,E)=\C.
$$

An ordered set
$E_0,\dots \dots,E_n$ of exceptional objects  is said to
be {\it exceptional} if for any~$i$
$$
\Hom^i(E_j,E_k)=0\quad\text{when}\quad j>k.
$$
An exceptional collection is
{\it full} if it generates the derived category.

One may choose to work with more specific classes of varieties:

\pph{\bf Fano varieties: cellular, minimal, Tate.}
An $m$-dimensional Fano variety $F$ is said to be \emph{minimal}
 if its cohomology
is as small as it can be ($H^{2k+1}(F,\Z)=0,H^{2k}(F,\Z)=\Z$).
A Fano  is said to be \emph{Tate}
if its motive has no non-Tate constituents. A Fano $F$ is said to be \emph{cellular} if
$F$ is a union of affine spaces:
$F=\bigcup {\mathbb{A}}_j^{i(j)}, \;{\mathbb{A}}_{j_1} \bigcap{\mathbb{A}}_{j_2}
= \varnothing \text{ if } j_1 \ne j_2.$ Cellular Fanos and minimal Fanos are Tate.

It has been noted that Tate Fanos `tend to' possess full exceptional collections,
though no precise conjecture has been made, to our knowledge. It is believed
that Fanos with exceptional collections are Tate. It is also believed that a minimal
Fano should have a full exceptional collection.

Minimal Fano threefolds are known to have
full exceptional collections, by Beilinson, Kapranov,
Orlov and Kuznetsov.

\pph \defi We call a basis $v_0,\dots,v_n$
of a linear space $V$ endowed with a bilinear form $\chi$ {\it
semiorthonormal} if $\chi(v_i,v_j)=0$ when $i>j$, and
$\chi(v_k,v_k)=1$ for all~$k$.

The classes of the elements of a full exceptional collection in
$K_0(F) \otimes \Q$ form a semiorthonormal basis with respect
to the Riemann-Roch form $\chi ([O_1],[O_2])=\sum (-1)^i \dimens \Hom ^i (O_1,O_2)$.

Nogin described semiorthogonal bases for minimal Fano threefolds in \cite{No}.

\bigskip

The statements of
\ref{par-1}, \ref{par-2}, though not necessary for the proof,
are worth to be kept in mind.

\pph \label{par-1} \bf The Coxeter element. \rm Let $A$ be the matrix of the
bilinear form $\chi$ in a semiorthonormal basis $v_0,\dots,v_n$
of a space $V$.

\bf (s) \rm Let  $\chi_s$ be the symmetrization of  $\chi$, that is
$ \chi_s(w_1,w_2)=\chi(w_1,w_2)+\chi(w_2,w_1),$
and let $I_0,\dots,I_n$ be the reflections with respect to the
vectors  $v_0,\dots,v_n$ in the orthogonal space $(V,\chi_s)$,
$
I_j\:v \mapsto v-\chi_s(v,v_j)v_j.
$
Then one has
$
I_0I_1\dots I_n=-A^{-1}A^t.
$

\bf (a) \rm Let  $\chi_a$ be the alternation of $\chi$, that
is
$
\chi_a(w_1,w_2)=\chi(w_1,w_2)-\chi(w_2,w_1),
$
and let $I_0,\dots,I_n$ be the reflections with respect to the
vectors $v_0,\dots,v_n$ in the symplectic space $(V,\chi_a)$.
Then one has
$
I_0I_1\dots I_n=A^{-1}A^t.
$

\qed

Now note that in the case  $V=K_0(F) \otimes \Q$  Serre's duality
 yields
$
\chi([O_1],[O_2])=\chi([O_2],[O_1 \otimes K_F[m]])
$for any pair of objects $O_1,O_2 \in \Ob \,D^b(F)$.

Let $A$ be the matrix of the form  $\chi$ in a basis that
consists of the classes of the elements of an exceptional
collection. Then the identity
$
x^tAy=y^tA(A^{-1}A^t)x
$
shows that in this basis the matrix of the operator defined on
 $K_0 \otimes \Q$ by the functor
$\otimes K_F [m]$ is  $A^{-1}A^t$. Hence one has

\pph \label{par-2} \propo In $K_0(F) \otimes \Q$:

 for $F$ odd-dimensional,
the product of all orthogonal reflections with respect to the
classes of the elements of an exceptional collection is  the
multiplication by~$[K_F]$; for $F$ even-dimensional,
 the product of all symplectic reflections with respect
to the classes of the elements of an exceptional collection is
the multiplication by~$[K_F]$;

\pph {\bf  The conjecture.}  Let $F$ be an
odd-dimensional (resp. even-dimensional) variety  $F$ with a full
exceptional collection  $\mathbf{E}=\langle E_1,E_2,\dots,E_n \rangle$.
Let $\chi^\bullet$~be the symmetrized (resp. alternated)
Riemann-Roch form. Consider the linear space $\overline{V}=(K_0(F)\otimes
\Q)/ \Ker \chi^\bullet$. The respective non-degenerate form
will also be denoted by
$\chi^\bullet$.

Let $\mathbf{x}=\langle x_1,\dots,x_n \rangle$ be an ordered set of
different points on $\A1$. Put
$U=\A1\setminus\{x_1,\dots,x_n\}$ and consider the
representation
$$
\phi_{\mathbf{E,\bdx}} :\pi_1(U^{an}) \to O(\overline{V},\chi^\bullet)
\quad\text{(respectively,}\quad \pi_1(U^{an})\to Sp(\overline{V},\chi^\bullet)),
$$
determined by requiring that the loop around the $i$-th point act
by reflection with respect to $[E_i]$. This representation
defines a local system  ${\cal L}_{\bdE,\bdx}$ on~$U$.

The conjecture that was stated by Dubrovin
in \cite{Du} (see also preceding discussion in \cite{Za})
says,
roughly, that if $\bdx$ is chosen
to be an ordered set of the critical values of
the so-called Landau--Ginsburg potential $u$, then ${\cal L}_{\bdE,\bdx}$
is the monodromy that arises in the local system of the middle cohomology
of the level set of the potential, $R^{mid}u_*(\Z)$.

As the notion of a Landau--Ginsburg model has no rigorous definition yet,
one interprets this statement by relating the local system to the monodromy
of Dubrovin's second structural connection.

\section{Regularized quantum $D$--module.}

As above, let $F$ be an $m$--dimensional
Fano variety of index $d$, so that
$-K_F=dH$.
Consider the matrix $M_{-K_F}$ of quantum multiplication by $-K_F$.
It has entries
in $\Q[q_i,q_i^{-1}]$, where $q_i$'s correspond, as usual, to
the generators of the lattice of numerical classes of curves on $F$.
Let $h_i$ be the anticanonical
degree of the class $q_i$.
One may specialize the matrix to $M$ in $\mathrm{Mat }\,(\Q [t,t^{-1}])$
by sending $q_i$  to $t^{h_i}$.
There is no need to do that when
$H^2(F)=\Z$ which we will assume henceforth. The anticanonical quantum $D$--module
on $\Gm$ is given by the connection
$t\frac{\partial}{\partial t} \eta = \eta M$, and the regularized quantum
$D$--module is the Fourier
transform of its [middle] extension to $\A1$.

The following is a
version of the
conjecture.

\pph  We shall say that the \bf exceptional collection/vanishing cycles
conjecture
holds for $F$ \rm if:

i) the regularized quantum $D$--module is regular singular (this is expected to hold
for any $F$), and the finite singularities $x_1, \dots x_n$ are simple;

ii) $F$ has an exceptional collection $\bdE$, and  $x_1, \dots x_n$
can be ordered to a tuple $\bdx$ so that
the local system $\L_{\bdE, \bdx}$ is isomorphic to the monodromy
of the regularized quantum $D$--module.

\pph Regularized quantum $D$--modules of the minimal Fano threefolds
$\P ^3, Q, V_5,V_{22}$ have been
identified with $d$--pullbacks of Pickard--Fuchs equations in
twisted Kuga--Sato families over the curves $X_0(N)/W_N$ where
$N=2,3,5,11$ respectively (${N=-K_F^3/2d^2}$,
see details in \cite{Go-CPMD}). In particular,
the sets $\mathbf{x}$ for the minimal Fano threefolds $F$ above have
cardinality $4$ and are
formed by the $d$--th roots of the elliptic points  $X_0(N)/W_N$, see below.

\section{$V_{22}$}

\pph \theo
The variety [rather, family of varieties] $V_{22}$ satisfies the
exceptional collection/vanishing
cycles conjecture.

\proof Kuznetsov showed in \cite{Ku96} that $V_{22}$ possesses
exceptional collections. An instance is $\mathbf{E}=({\mathcal O},S^*,E^*,\Lambda^2S^*)$, in
Kuznetsov's notation. The matrices of $\chi, \chi^\bullet=\chi^{sym}$
in this basis have the
form
$$
X=\left(\begin{array}{rrrr}
1 & 7 & 8 & 18 \\
0 & 1 & 4 & 13 \\
0 & 0 & 1 & 4  \\
0 & 0 & 0 & 1
\end{array}\right),\:  X^{sym}=\left(\begin{array}{rrrr}
2 & 7 & 8 & 18 \\
7 & 2 & 4 & 13 \\
8 & 4 & 2 & 4  \\
18 & 13 & 4 & 2
\end{array}\right).
$$

On the other hand, the monodromy of the regularized quantum $D$--
module was described in \cite{Go-CPMD} as follows.

\pph \label{identification}
Let $X_0(N)^\circ$ stand for $X_0(N)-\{
\text{cusps} \} -\{ \text{elliptic points} \}$.
Let $\phi$ be the
tautological projective representation
$\phi \colon
\pi _1(X_0(N)^\circ) \longrightarrow PSL_2(\Z).$ The monodromy that acts on $H^1$
of the fiber of `a universal
elliptic curve' is given by a lift of $\phi$ to a linear
representation
$$\bar \phi \colon \gamma  \mapto \begin{pmatrix}
   a & b \\
   c & d
\end{pmatrix}, \; c=0 \mod N.$$
In a suitable basis its symmetric square is
$$
\psi\colon \gamma \mapto Sym_N^2 \, \phi (\gamma)=\begin{pmatrix}
   d^2 & 2cd & -c^2/N \\
   bd & bc+ad & -ac/N \\
   -Nb^2 & -2Nab &a^2
\end{pmatrix}.
$$

Let $W$ be the Atkin--Lehner involution given by the action
of $W=\begin{pmatrix}
0 & 1 \\
-N & 0
\end{pmatrix}$ on $X_0(N)$.
Delete the $W$-invariant points from $X_0(N)^\circ$ and
let ${X_0(N)^W}^\circ$ be the quotient of the resulting curve by
$W$.
The fundamental group of ${X_0(N)^W}^\circ$ is  then generated by $ \pi
_1(X_0(N)^\circ)$ and a loop $\iota$ around the point that is the
image of a point $s$ on the upper halfplane stabilized by $W$.
Extend $\psi$ to $\iota$, by setting $\psi (\iota)=I$
with $I=\begin{pmatrix}
   0 & 0 & 1 \\
   0 & 1 & 0 \\
   1 & 0 & 0
\end{pmatrix}.$
For a Fano $F$ with $\displaystyle{d=1, \;N=\frac{-K^3_F}{2d^2}}$,
the resulting representation is the monodromy representation of the regularized
quantum $D$--module.

\bigskip
\bigskip

Having said all that, and made choices, one may describe the
monodromy of the regularized quantum D-module for $V_{22}$ as
follows. Let
$\gamma_{12}, \gamma_{13},\gamma_{14},\gamma_{23}, \gamma_{24},\gamma_{34} \in \Gamma_0(11)$
respectively
denote
$$\left (\begin {array}{cc} 4&1\\\noalign{\medskip}11&3\end {array}
\right ),\left (\begin {array}{cc} 6&1\\\noalign{\medskip}11&2
\end {array}\right ),\left (\begin {array}{cc} 15&2
\\\noalign{\medskip}22&3\end {array}\right ),\left (\begin {array}
{cc} 7&1\\\noalign{\medskip}-22&-3\end {array}\right ),\left (
\begin {array}{cc} 23&3\\\noalign{\medskip}-77&-10\end {array}\right )
,\left (\begin {array}{cc} 8&1\\\noalign{\medskip}-33&-4\end {array}
\right )$$
so that $\gamma_{12}\gamma_{23}=\gamma_{13},
\gamma_{12}\gamma_{24}=\gamma_{14},
\gamma_{23}\gamma_{34}=\gamma_{24}$.
Notice that $Tr \, \gamma_{ij} = X_{ij}$.

Let $W=\left (\begin {array}{cc} 0&-1\\\noalign{\medskip}11&0\end
{array} \right )$ as above. Then the images of the points on the upper half
plane stabilized by $W,W\gamma_{12},W\gamma_{13},W\gamma_{14}$ are
the four elliptic points on $X_0(11)^W$. The monodromies around
these points are, respectively, the reflections $ I, I \psi
(\gamma_{12}),I \psi (\gamma_{13}),I \psi (\gamma_{14})$. These
monodromies are orthogonal with respect to the form $U=\left (\begin
{array}{ccc} 0&0&-1\\\noalign{\medskip}0&-22&0
\\\noalign{\medskip}-1&0&0\end {array}\right )$, that is,
$I^tUI=\psi(\gamma_{ij})^tU\psi(\gamma_{ij})=U$. The length $2$
vectors of the reflections are, respectively,
$v_1=(-1,0,1),v_2=(-4,1,3),v_3=(-6,1,2),v_4=(-15,2,3)$. One notices
that $(v_i,v_j)={X^{sym}}_{ij}$. Therefore, the exceptional collection/vanishing
cycles
conjecture holds for $V_{22}$.

\section{$V_5$}

The proof in the case of the variety $V_5$ goes along the same lines,
except that the index is now $2$, and the monodromy of the regularized
quantum D--module is realized as the restriction of $\psi$ to an index $2$
subgroup of $\pi_1({X_0(5)^W}^\circ)$ which corresponds to its double cover
ramified over the unique cusp and the elliptic point that is the image of the
two elliptic points on  $X_0(5)$.

Let
$\gamma_{12}, \gamma_{13},\gamma_{14},\gamma_{23}, \gamma_{24},\gamma_{34} \in \Gamma_0(5)$
respectively
denote
$$
\left (\begin {array}{cc} 2&1\\\noalign{\medskip}5&3\end {array}
\right ),\left (\begin {array}{cc} 3&1\\\noalign{\medskip}5&2
\end {array}\right ),\left (\begin {array}{cc} 6&1\\\noalign{\medskip}
5&1\end {array}\right ),\left (\begin {array}{cc} 4&1
\\\noalign{\medskip}-5&-1\end {array}\right ),\left (\begin {array}
{cc} 13&2\\\noalign{\medskip}-20&-3\end {array}\right ),\left (
\begin {array}{cc} 7&1\\\noalign{\medskip}-15&-2\end {array}\right ).
$$
They generate an index $2$ subgroup of $\Gamma_0(5)$; the respective cover ramifies
over the $2$ cusps and the $2$ elliptic points on $X_0(5)$.

Let $W=\left (\begin {array}{cc} 0&-1\\\noalign{\medskip}5&0\end
{array} \right )$. Then the images of the points on the upper half
plane stabilized by $W,W\gamma_{12},W\gamma_{13},W\gamma_{14}$ are
the four elliptic points on the double cover of $X_0(5)^W$.
The monodromies around
these points are, respectively, the reflections $ I, I \psi
(\gamma_{12}),I \psi (\gamma_{13}),I \psi (\gamma_{14})$. These
monodromies are orthogonal with respect to the form $U=\left (\begin
{array}{ccc} 0&0&-1\\\noalign{\medskip}0&-10&0
\\\noalign{\medskip}-1&0&0\end {array}\right )$. The length $2$
vectors of the reflections are, respectively,
$v_1=(-1,0,1),v_2=(-2,1,3),v_3=(-3,1,2),v_4=(-6,1,1)$. Then
 $(v_i,v_j)={X^{sym}}_{ij}$
with
$$
X=\left (\begin {array}{cccc} 1&5&5&7\\\noalign{\medskip}0&1&3&10
\\\noalign{\medskip}0&0&1&5\\\noalign{\medskip}0&0&0&1\end {array}
\right )
$$

Orlov proved that $\mathbf{E}= (\OO, Q, S^*, \OO (1))$ as in \cite{Or} is an exceptional
collection with $\chi ([E_i],[E_j])=X_{ij}$ as above.
Therefore, the exceptional collection/vanishing
cycles
conjecture holds for $V_{5}$.

\section{The quadric}

The index is now $3$, and the monodromy of the regularized
quantum D--module is realized as the restriction of $\psi$ to an index $3$
subgroup of $\pi_1({X_0(3)^W}^\circ)$ which corresponds to its cyclic $3$--cover
ramified over the unique cusp and the elliptic point that is the image of the unique
elliptic point on  $X_0(3)$.

Let
$\gamma_{12}, \gamma_{13},\gamma_{14},\gamma_{23}, \gamma_{24},\gamma_{34} \in \Gamma_0(3)$
respectively
denote
$$
\left (\begin {array}{cc} 2&1\\\noalign{\medskip}3&2\end {array}
\right ),\left (\begin {array}{cc} 4&1\\\noalign{\medskip}3&1
\end {array}\right ),\left (\begin {array}{cc} 13&2
\\\noalign{\medskip}6&1\end {array}\right ),\left (\begin {array}{cc}
5&1\\\noalign{\medskip}-6&-1\end {array}\right ),\left (\begin {array}
{cc} 20&3\\\noalign{\medskip}-27&-4\end {array}\right ),\left (
\begin {array}{cc} 7&1\\\noalign{\medskip}-15&-2\end {array}\right ).
$$
They generate an index $3$ subgroup of $\Gamma_0(3)$; the respective cover ramifies
over the $2$ cusps and the  elliptic point on $X_0(3)$.

Let $W=\left (\begin {array}{cc} 0&-1\\\noalign{\medskip}3&0\end
{array} \right )$ as above. Then the images of the points on the upper half
plane stabilized by $W,W\gamma_{12},W\gamma_{13},W\gamma_{14}$ are
the four elliptic points on the $3$--cover of $X_0(3)^W$.
The monodromies around
these points are, respectively, the reflections $ I, I \psi
(\gamma_{12}),I \psi (\gamma_{13}),I \psi (\gamma_{14})$. These
monodromies are orthogonal with respect to the form $U=\left (\begin
{array}{ccc} 0&0&-1\\\noalign{\medskip}0&-6&0
\\\noalign{\medskip}-1&0&0\end {array}\right )$. The length $2$
vectors of the reflections are, respectively,
$v_1=(-1, 0, 1),v_2=(-2, 1, 2),v_3=(-4, 1, 1),v_4=(-13, 2, 1).$ Then
 $(v_i,v_j)={X^{sym}}_{ij}$
with
$$
X=\left (\begin {array}{cccc} 1&4&5&14\\\noalign{\medskip}0&1&4&16
\\\noalign{\medskip}0&0&1&5\\\noalign{\medskip}0&0&0&1\end {array}
\right )
$$

Let $S$ be the spinor bundle on the quadric.
The collection $\mathbf{E}= ( \OO, S^*,\OO (1), \OO(2) )$ is an exceptional
collection  with $\chi ([E_i],[E_j])=X_{ij}$ as above (e.g. \cite{No}).
Thus, we have shown that the exceptional collection/vanishing
cycles conjecture holds for the $3$--dimensional quadric.

\section{$\P ^3$}

The index is now $4$, and the monodromy of the regularized
quantum D--module is realized as the restriction of $\psi$ to an index $4$
subgroup of $\pi_1({X_0(2)^W}^\circ)$ which corresponds to its cyclic $4$--cover
ramified over the unique cusp and the elliptic point that is the image of the unique
elliptic point on  $X_0(2)$.

Let
$\gamma_{12}, \gamma_{13},\gamma_{14},\gamma_{23}, \gamma_{24},\gamma_{34} \in \Gamma_0(2)$
respectively
denote
$$
\left (\begin {array}{cc} 3&1\\\noalign{\medskip}2&1\end {array}
\right ),\left (\begin {array}{cc} 9&2\\\noalign{\medskip}4&1
\end {array}\right ),\left (\begin {array}{cc} 19&3
\\\noalign{\medskip}6&1\end {array}\right ),\left (\begin {array}{cc}
5&1\\\noalign{\medskip}-6&-1\end {array}\right ),\left (\begin {array}
{cc} 13&2\\\noalign{\medskip}-20&-3\end {array}\right ),\left (
\begin {array}{cc} 7&1\\\noalign{\medskip}-22&-3\end {array}\right )
.
$$
They generate an index $4$ subgroup of $\Gamma_0(2)$; the respective cover ramifies
over  the $2$ cusps and the  elliptic point on $X_0(2)$.

Let $W=\left (\begin {array}{cc} 0&-1\\\noalign{\medskip}2&0\end
{array} \right )$. Then the images of the points on the upper half
plane stabilized by $W,W\gamma_{12},W\gamma_{13},W\gamma_{14}$ are
the four elliptic points on the $4$--cover of $X_0(2)^W$.
The monodromies around
these points are, respectively, the reflections $ I, I \psi
(\gamma_{12}),I \psi (\gamma_{13}),I \psi (\gamma_{14})$. These
monodromies are orthogonal with respect to the form $U=\left (\begin
{array}{ccc} 0&0&-1\\\noalign{\medskip}0&-4&0
\\\noalign{\medskip}-1&0&0\end {array}\right )$. The length $2$
vectors of the reflections are, respectively,
$v_1=(-1, 0, 1),v_2=(-3, 1, 1),v_3=(-9, 2, 1),v_4=(-19, 3, 1)$. Then
 $(v_i,v_j)={X^{sym}}_{ij}$
with
$$
X=\left (\begin {array}{cccc} 1&4&10&20\\\noalign{\medskip}0&1&4&10
\\\noalign{\medskip}0&0&1&4\\\noalign{\medskip}0&0&0&1\end {array}
\right )
$$

By the Kodaira vanishing theorem,
$\mathbf{E}= ( \OO, \OO (1),\OO (2), \OO(3) )$ is an exceptional
collection with $\chi ([E_i],[E_j])=X_{ij}$ as above.
This is another explanation why the exceptional collection/vanishing
cycles conjecture holds for $\P ^3$.

\pph \bf Remarks and acknowledgements. \rm A generating
set of Gromov--Witten invariants for $V_{22}$ had been computed by Kuznetsov
(published in \cite{BaMa}),
and a set for $V_5$, by Beauville \cite{Bea}, thus providing the input for constructing
Dubrovin's connection. The truth of the conjecture for projective spaces and
odd--dimensional quadrics is well--known and is probably due to Kontsevich. See
\cite{Gu} for a version with the first structural
connection and the Stokes matrix, or
an $l$--adic version in \cite{Go-RRV}.
An issue we do not even mention, that of constructibility at the categorial level,
is worked out by Polishchuk in \cite{Polishchuk-2007}.

I thank Boris Dubrovin, Tamara Grava,
Yuri Manin and Dmitri Orlov for discussions of the subject, and Constantin Shramov
for proofreading this note.

\end{document}